\documentclass[twocolumn]{autart}

\usepackage{dsfont}
\usepackage{amsfonts}
\usepackage{mathrsfs}
\usepackage{}
\usepackage{colortbl}
\usepackage{amssymb,amsmath,graphicx}
\usepackage{float}
\usepackage{multirow}

\newtheorem{Remark}{Remark}[section]

\newtheorem{Problem}{Problem}
\newenvironment{Proof}{\noindent{\em Proof:\/}}{\hfill $\Box$\par}
\newtheorem{Theorem}{Theorem}[section]
\newtheorem{Proposition}{Proposition}[section]
\newtheorem{Lemma}{Lemma}[section]

\newtheorem{Assumption}{Assumption}

\newcommand{\EQQ}{\begin{eqnarray*}}
\newcommand{\ENN}{\end{eqnarray*}}
\newcommand{\EQ}{\begin{eqnarray}}
\newcommand{\EN}{\end{eqnarray}}

\begin{document}

\begin{frontmatter}

\title{Event-Triggered Cooperative Robust Practical Output Regulation for a Class of Linear  Multi-Agent Systems
\thanksref{footnoteinfo}} 

\thanks[footnoteinfo]{This work has been supported in part by the Research Grants
Council of the Hong Kong Special Administration Region under grant
No. 14200515, and in part by National Natural Science Foundation of
China under Project 61633007.
Corresponding author: Jie Huang.
}

\author[Liu-Huang]{Wei Liu}\ead{wliu@mae.cuhk.edu.hk},
\author[Liu-Huang]{Jie Huang}\ead{jhuang@mae.cuhk.edu.hk}

\address[Liu-Huang]{Department of Mechanical and Automation Engineering, The Chinese University of Hong Kong, Shatin, N.T., Hong Kong}


\begin{abstract}
In this paper, we consider the event-triggered cooperative robust practical output regulation problem for a class of linear  minimum-phase multi-agent systems. We first convert our problem into the cooperative robust practical stabilization problem of a well defined augmented system
Based on  the distributed internal model approach. Then, we design a distributed event-triggered output feedback control law together with a distributed output-based  event-triggered mechanism to stabilize the augmented system, which leads to the solvability of the cooperative robust practical output regulation problem of the original plant. Our distributed control law can be directly implemented in a digital platform provided that the distributed triggering mechanism can monitor the continuous-time output information from neighboring agents.
Finally, we illustrate  our design by an example.
\end{abstract}

\end{frontmatter}

\section{Introduction}

Over the past few years, event-triggered control has attracted extensive attention from the control community. Compared with conventional periodical time-triggered control, event-triggered control samples the system's state or output aperiodically and  reduces the number of control task executions while maintaining the control performance.   As noted in \cite{Heemels1},   event-triggered control is reactive and generates sensor sampling and control actuation when, for instance, the plant state or output deviates more than a certain threshold from a desired value.
 Various event-triggered control problems have been studied for several types of systems in \cite{Donkers2012,Girard1,LiuT1,LiuHuang2017,Tabuada1,Tallapragada1} and the references therein. One of the main challenges for  event-triggered control is to avoid the Zeno behavior, which means that the execution times become arbitrarily close and result in an accumulation point \cite{Tabuada1}. In \cite{Tabuada1},   the  stabilization problem of a given system by a state-based event-triggered control law was linked to
 the input-to-state stabilizability (ISS) of the system, and it was shown that the Zeno behavior can be excluded when the measurement state is not equal to zero.
 In \cite{Girard1}, a dynamic state-based  event-triggered mechanism was further proposed to study the stabilization problem for the same class of nonlinear systems as that in \cite{Tabuada1}.
Reference \cite{Tallapragada1} studied the asymptotic tracking problem of a control system by a stated-based event-triggered controller and it gave
a stated-based event-triggered controller that was able to achieve the uniform ultimate boundedness of the tracking error.
 Reference \cite{LiuT1} further studied the robust stabilization problem for a class of systems subject to external disturbances by a state-based event-triggered control based on the small-gain theorem.
  In \cite{Donkers2012}, an output-based event-triggered control law was proposed to analyze the stability and $\mathcal{L}_{\infty}$ performance for a class of linear systems.
  In \cite{LiuHuang2017},  the robust practical output regulation problem for a class of linear systems was studied by an output-based event-triggered control law.
  Some other contributions can also be found in \cite{Abdelrahim1,Dolk2016,Postoyan1} etc.

In this paper, we further consider the cooperative robust practical output regulation problem for a class of  linear minimum-phase multi-agent systems by a dynamic distributed output-based event-triggered control law. 
The problem can be viewed as an extension of the result in \cite{LiuHuang2017} from a single system to a multi-agent system.
Compared with \cite{LiuHuang2017},   the main challenge for this paper is that the control of each subsystem and hence the triggering mechanism are subject to some communication constraints described by a digraph. We need to specifically design an output-based event-triggered control law that satisfies the communication constraints. Thus, our
 problem is  more challenging than that of \cite{LiuHuang2017}. In what follows, a control law and an event-triggered mechanism satisfying the communication constraints are called a distributed control law and a distributed event-triggered mechanism, respectively.


Several other cooperative control problems of multi-agent systems by event-triggered control have also been studied. For example,
the consensus problem of single-integrator and double-integrator multi-agent systems by event-triggered control was studied in \cite{Dimarogonas1,Fan1,Seyboth1}, and  \cite{Li1,Mu1}, respectively. Reference \cite{Zhang1} studied the consensus problem for the general linear multi-agent systems by an observer-based output feedback event-triggered distributed control law. Reference \cite{Wang4} studied the cooperative output regulation problem for linear multi-agent system by a  centralized event-triggered mechanism. Reference \cite{HuW2} further studied the cooperative output regulation problem for linear multi-agent systems by a distributed event-triggered mechanism based on the feedforward design method.
Other relevant results for the general linear multi-agent systems can be found in \cite{Cheng1,Wang5,Zhu1} etc.

Compared with the existing results on event-triggered cooperative control problems,  this paper needs to overcome some specific challenges. First,  our problem formulation generalizes the existing event-triggered  cooperative control problems, say, in \cite{Cheng1,Dimarogonas1,Fan1,Li1,Mu1,Seyboth1,Zhang1,Zhu1}, in the sense that we  achieve not only asymptotical tracking but also disturbance rejection. Second, since our system contains  unknown parameters, we need to adopt the distributed internal model approach, which leads to a robust stabilization problem
 for a much more complicated augmented system.
 Third,  our event-triggered mechanism is output-based and distributed in the sense that the event-triggered mechanism of each agent only depends on the output information of its neighbors and itself. Finally,  our control law is piecewise constant, which lends itself to a direct  implementation in a digital platform.



Throughout this paper, we use the following notation: $\mathbb{Z}^{+}$ denotes the set of all nonnegative integers. For any column vectors $a_i$, $i=1,...,s$, we denote
$\mbox{col}(a_1,...,a_s)=[a_1^T,...,a_s^T]^T$. For any matrices $X\in\mathbb{R}^{n\times m}$, we denote $\mbox{vec}(X)=[X_{1}^{T},\cdots,X_{m}^{T}]^{T}$ where $X_{i}$ with $i=1,\cdots,m$ is the $i$th column of $X$.
The notation $\|x\|$ denotes the Euclidean norm of vector $x$. The notation $\|A\|$ denotes the induced norm of matrix $A$ by the Euclidean norm. The notation $\otimes$ denotes the Kronecker product of matrices. Denote the base of the natural logarithm  by $\mathbf{e}$.
Denote the maximum eigenvalue and the minimum eigenvalue of a symmetric real matrix $A$ by $\lambda_{\max}(A)$ and $\lambda_{\min}(A)$, respectively.

\section{Problem formulation and preliminaries}\label{PF}
Consider the following linear multi-agent systems
\begin{equation}\label{system1}
\begin{split}
&\dot{z}_{i}=A_{1i}(w)z_{i}+A_{2i}(w)\xi_{1i}+E_{0i}(w)v\\
&\dot{\xi}_{si}=\xi_{(s+1)i}, \ \ s=1,\cdots,r-1\\
&\dot{\xi}_{ri}\!=\!A_{3i}(w)z_{i}\!+\!\sum_{s=1}^{r}c_{si}(w)\xi_{si}\!+\!E_{ri}(w)v\!+\!b_{i}(w)u_{i}\\
&y_{i}=\xi_{1i},\ i=1,\cdots,N,
\end{split}
\end{equation}
where $z_{i}\in\mathbb{R}^{n_{i}-r}$ and $\xi_{i}=\mbox{col}(\xi_{1i},\cdots,\xi_{ri})\in\mathbb{R}^{r}$ are the states, $u_{i}\in \mathbb{R}$ is the input, $y_{i}\in \mathbb{R}$ is the output,
$w\in\mathbb{R}^{n_{w}}$ is an unknown parameter vector, $b_{i}(w)>0$ for all $w\in\mathbb{R}^{n_w}$, and $v(t)\in\mathbb{R}^{n_{v}}$ is an exogenous signal representing both reference input and disturbance and is assumed to be generated by the following linear system
\begin{equation}\label{exosystem1}
\begin{split}
  \dot{v} = Sv
\end{split}
\end{equation}
where $S$ is some known constant matrix.  Let $y_{0}=F(w)v\in\mathbb{R}$ be the output of the exosystem \eqref{exosystem1}. Then, the regulated error of each subsystem is defined
 as  $e_{i}=y_{i}-y_{0}$ for $i=1,\cdots,N$.

For $i=1,\cdots,N$ and $s=1,\cdots,r$, let
$A_{1i}(w)=A_{1}+w_{A_{1i}}\in\mathbb{R}^{(n_{i}-r)\times (n_{i}-r)}$, $A_{2i}(w)=A_{2}+w_{A_{2i}}\in\mathbb{R}^{(n_{i}-r)\times 1}$,
$A_{3i}(w)=A_{3}+w_{A_{3i}}\in\mathbb{R}^{1\times (n_{i}-r)}$, $E_{0i}(w)=E_{0}+w_{E_{0i}}\in\mathbb{R}^{(n_{i}-r)\times n_{v}}$,
$E_{ri}(w)=E_{r}+w_{E_{ri}}\in\mathbb{R}^{1\times n_{v}}$, $F(w)\!=\!F\!+\!w_{F}\in\mathbb{R}^{1\times n_{v}}$,
$c_{si}(w)=c_{s}\!+\!w_{c_{si}}\in\mathbb{R}$, $b_{i}(w)\!=\!b+w_{b_{i}}\in\mathbb{R}$, $w_{A_{1}}\!\!=\!\mbox{col}(\mbox{vec}(w_{A_{11}}),\cdots,\mbox{vec}(w_{A_{1N}}))$,
$w_{A_{2}}\!\!=$ $\mbox{col}(\mbox{vec}(w_{A_{21}}),\cdots\!,\mbox{vec}(w_{A_{2N}}))$, $w_{A_{3}}\!\!=\!\mbox{col}(\mbox{vec}(w_{A_{31}}),\cdots$, $\mbox{vec}(w_{A_{3N}}))$,
$w_{E_{0}}\!\!=\!$ $\mbox{col}(\mbox{vec}(w_{E_{01}}),\cdots,\mbox{vec}(w_{E_{0N}}))$, $w_{E_{r}}$ $=\mbox{col}(\mbox{vec}(w_{E_{r1}}),\cdots,\mbox{vec}(w_{E_{rN}}))$, $w_{c_{s}}= \mbox{col}(w_{c_{s1}},\cdots$, $w_{c_{sN}})$, $w_{b}=\mbox{col}(w_{b_{1}},\cdots,w_{b_{N}})$,
$w=\mbox{col}(w_{A_{1}},w_{A_{2}},w_{A_{3}}$, $w_{E_{0}},w_{E_{r}},w_{c_{1}},\cdots,w_{c_{r}},w_{b},w_{F}^{T})$,
where $A_{1},A_{2},A_{3},E_{0}$, $E_{r},F,c_{s}$ and $b$  denote the nominal values, and $w$ denotes the unknown parameter.

System (\ref{system1}) is called a linear multi-agent system in the normal form and it is said to be minimum phase if the matrix $A_{1i}(w)$ is Hurwitz for all $w$. Systems (\ref{system1}) and (\ref{exosystem1}) together can be viewed as a multi-agent system with (\ref{exosystem1}) as the leader system, and the $N$ subsystems of (\ref{system1}) as $N$ followers.
It is noted that the cooperative robust output regulation problem for the multi-agent system composed of \eqref{system1} and \eqref{exosystem1} was studied by the continuous distributed control law  in \cite{Su6}.

As in \cite{Su6},  given the plant (\ref{system1}) and  the exosystem (\ref{exosystem1}),  we can define a digraph
$\bar{\mathcal{G}}=(\bar{\mathcal{V}},\bar{\mathcal{E}})$\footnote{See \cite{Su6}  for a summary of digraph.},
where $\bar{\mathcal{V}}=\{0,1\cdots,N\}$ with $0$ associated with
the leader system and with $ i = 1,\cdots,N$ associated with the $N$
followers, respectively, and $\bar{\mathcal{E}}\subseteq
\bar{\mathcal{V}}\times\bar{\mathcal{V}}$ for all $t\geq0$. For all
$t \geq 0$, each $j=0,1,\cdots,N$, $i=1,\cdots,N$, and $i \neq j$, $(j,i)\in\bar{\mathcal{E}} $  if and only if the control
$u_{i} (t)$ can make use of $y_{i} (t) -y_{j} (t)$ for feedback
control.
Let $\bar{\mathcal{N}}_{i}=\{j, (j,i)\in\bar{\mathcal{E}}\}$ denote the neighbor set of node $i$.

To define our control law, we recall that the adjacency matrix of the digraph
$\bar{\mathcal{G}}$ is a nonnegative matrix
$\bar{\mathcal{A}} =[\bar{a}_{ij} ]\in \mathbb{R}^{(N+1)\times
(N+1)}$ where $\bar{a}_{ii} =0$, $\bar{a}_{ij} =1\Leftrightarrow
(j,i)\in\mathcal{\bar{E}}$, and $\bar{a}_{ij} =0\Leftrightarrow (j,i)\notin\mathcal{\bar{E}}$ for $i,j=0,1,\cdots,N$.
Like in \cite{Su6},  define the virtual output for the $i$th subsystem as $e_{vi}(t)=\sum_{j=0}^{N}\bar{a}_{ij}(y_{i}(t)-y_{j}(t))$ for $i=1,\cdots,N$.
Let $e_{0}=0$, $e=\mbox{col} (e_{1},\cdots,e_{N})$, $e_{v}=\mbox{col} (e_{v1},\cdots,e_{vN})$, and
$H=[h_{ij}]_{i,j=1}^{N}$ with $h_{ii} =\sum_{j=0}^{N}\bar{a}_{ij}$ and $h_{ij}=-\bar{a}_{ij}$ for $i\neq j$.
It can be verified that $e_{v}=He$.

For any $k\in\mathbb{Z}^{+}$  and $i=1, \cdots,N$, consider the following control law
\begin{equation}\label{ui1}
\begin{split}
  &u_{i}(t)=F_{1}\eta_{i}(t_{k}^{i})+F_{2}\zeta_{i}(t_{k}^{i})\\
  &\dot{\eta}_{i}(t)= G_{1}\eta_{i}(t)+G_{2}\eta_{i}(t_{k}^{i})+G_{3}\zeta_{i}(t_{k}^{i})\\
  &\dot{\zeta}_{i}(t)= G_{4}\zeta_{i}(t)+G_{5}e_{vi}(t_{k}^{i})\\
\end{split}
\end{equation}
where  $F_{1}$, $F_{2}$, $G_{1},\cdots, G_{5}$ are some real matrices with proper dimensions, the $\eta_{i}$ subsystem is the so-called internal model, the $\zeta_{i}$ subsystem is a dynamic compensator, and $t_{k}^{i}$ denotes the triggering time instants of agent $i$ with $t_{0}^{i}=0$ and is generated by  the following event-triggered mechanism:
\begin{equation}\label{trigger1}
\begin{split}
t_{k+1}^{i}\!\!=\!\inf\{t\!>\!t_{k}^{i}|h_{i}(\tilde{e}_{vi}(t),\tilde{\eta}_{i}(t),\tilde{\zeta}_{i}(t),e_{vi}(t),\zeta_{i}(t),t)\!\geq\!\delta\}
\end{split}
\end{equation}
where $h_{i}(\cdot)$ is some nonlinear function, $\delta>0$ is a constant and  $\tilde{e}_{vi}(t)=e_{vi}(t_{k}^{i})-e_{vi}(t)$, $\tilde{\eta}_{i}(t)=\eta_{i}(t_{k}^{i})-\eta_{i}(t)$, $\tilde{\zeta}_{i}(t)=\zeta_{i}(t_{k}^{i})-\zeta_{i}(t)$.
\begin{Remark}
The event-triggering mechanism \eqref{trigger1} and the control law \eqref{ui1} are both distributed and output-based, and we call \eqref{ui1} as the dynamic distributed output-based event-triggered control law.
\end{Remark}
Our problem is described as follows.
\begin{Problem}\label{LCORPS}
 Given the plant (\ref{system1}), the exosystem (\ref{exosystem1}), any compact sets $\mathbb{V}\in\mathbb{R}^{n_{v}}$ with $0\in\mathbb{V}$ and $\mathbb{W} \subset \mathbb{R}^{n_{w}} $ with $0\in\mathbb{W}$,  a  digraph $\bar{\mathcal{G}}$, and any $\epsilon>0$, design an event-triggered mechanism of the form \eqref{trigger1} 
  and a control law of the form (\ref{ui1})  such that the resulting closed-loop system has the following properties: for any $v\in\mathbb{V}$, $w\in\mathbb{W}$, and any initial conditions $z_{i}(0)$, $\xi_{i}(0)$, $\eta_{i}(0)$, $\zeta_{i}(0)$ for $i=1,\cdots,N$,
  \begin{description}
    \item[a)] the trajectory of the closed-loop system exists and is bounded for all time $t\geq0$;
    \item[b)] $\lim_{t \to +\infty}\|e(t)\|\leq\epsilon$.
  \end{description}
\end{Problem}
\begin{Remark}\label{RemarkProb1a}
We call Problem \ref{LCORPS} as the event-triggered cooperative robust practical output regulation problem. Clearly, Problem \ref{LCORPS} is more challenging than the cooperative robust output regulation problem studied in \cite{Su6} in at least three ways.
First, we not only need to design a control law, but also an
event-triggered
mechanism. Moreover,  both the control law and the event-triggered mechanism must be distributed in the sense that both the control law and
the event-triggered mechanism of each subsystem can only make use of the
information of its neighbors and itself.
Second,  since the closed-loop system is
hybrid under event-triggered control, the stability analysis of the closed-loop system is much more complicated than that of the  closed-loop system in \cite{Su6} which is smooth.
Third, we need to make a special effort to design an appropriate event-triggered mechanism to prevent the Zeno phenomenon.
It is also interesting to note that, unlike the control laws in  \cite{HuW2,Wang4,Wang5,Zhang1},
 the specific form of the control law \eqref{ui1} lends itself to a direct digital implementation as follows:
\begin{equation}\label{ui2}
\begin{split}
  u_{i}(t)\!=&F_{1}\eta_{i}(t_{k}^{i})+F_{2}\zeta_{i}(t_{k}^{i})\\
  \eta_{i}(t_{k+1}^{i})\!=&\mathbf{e}^{G_{1}(t_{k+1}^{i}-t_{k}^{i})}\eta_{i}(t_{k}^{i})+\\
  &(G_{2}\eta_{i}(t_{k}^{i})\!+\!G_{3}\zeta_{i}(t_{k}^{i}))\!\!\!\int_{t_{k}^{i}}^{t_{k+1}^{i}}\!\!\!\mathbf{e}^{G_{1}(t_{k+1}^{i}\!-\!\tau)}d\tau\\
  \zeta_{i}(t_{k+1}^{i})\!=&\mathbf{e}^{G_{4}(t_{k+1}^{i}-t_{k}^{i})}\zeta_{i}(t_{k}^{i})+\\
  &G_{5}e_{vi}(t_{k}^{i})\int_{t_{k}^{i}}^{t_{k+1}^{i}}\!\mathbf{e}^{G_{4}(t_{k+1}^{i}-\tau)}d\tau\\
\end{split}
\end{equation}
for $t\in[t_{k}^{i},t_{k+1}^{i})$ with $k\in\mathbb{Z}^{+}$ and $i=1,\cdots,N$.
\end{Remark}
\begin{Remark}\label{RemarkProb1b}
It is noted that the event-triggered mechanism \eqref{trigger1} relies on the continuous-time measurement output $e_{vi}(t)$, thus costing communication resources.
It is interesting to further consider how to design novel triggering mechanisms that avoid continuous communication by appropriate samplings.
\end{Remark}


To solve our problem,   we need some assumptions as follows.
\begin{Assumption}\label{Ass1}
The exosystem is neutrally stable, i.e., all the eigenvalues of S are semi-simple with zero real parts.
\end{Assumption}

\begin{Assumption}\label{Ass2}
For $i=1,\cdots,N$,  each subsystem of \eqref{system1} is minimum phase for any $w\in\mathbb{W}$, i.e. $A_{1i}(w)$ is Hurwitz for all $w\in\mathbb{W}$.
\end{Assumption}

\begin{Assumption}\label{Ass3}
The digraph $\bar{\mathcal{G}}$ is connected, i.e., it contains a directed spanning tree with node $0$ as the root.
\end{Assumption}

\begin{Remark}
Assumption \ref{Ass1} has also been  used in \cite{LiuHuang2017}. Under Assumption \ref{Ass1}, for any $v(0)\in\mathbb{V}_{0}$ with $\mathbb{V}_{0}$ being a compact set, there exists another compact set $\mathbb{V}$ such that $v(t)\in\mathbb{V}$ for all $t\geq0$.
Assumptions \ref{Ass2} and \ref{Ass3} are the same as Assumptions 6 and 7 of \cite{Su6}.
Under Assumption \ref{Ass3}, by Lemma 4 in \cite{Hu1}, $H$ is nonsingular and all eigenvalues of $H$ have positive real parts.
\end{Remark}

\section{The existing result}\label{ER}
In this section, we will summarize the result of \cite{Su6}, where the cooperative robust output regulation problem for the system \eqref{system1} has been studied by a continuous-time distributed output feedback control law.



As shown in \cite{Su6}, under Assumptions \ref{Ass1}-\ref{Ass2}, we can define a set of dynamic compensators
\begin{equation}\label{dotetai}
\begin{split}
  \dot{\eta}_{i}=M\eta_{i}+Qu_{i},\  i = 1, \cdots, N
\end{split}
\end{equation}
 where
$(M,Q)$ is any controllable pair with $M\in \mathbb{R}^{l\times l}$ being a Hurwitz matrix, $Q\in \mathbb{R}^{l\times 1}$ being a column vector and $l$ being the degree of the minimal polynomial of $S$. This set of dynamic compensators is called  the distributed internal model of the system (\ref{system1}).
Let
$\Phi=\begin{bmatrix}
        \begin{smallmatrix}
         0 & 1 & \cdots & 0 \\
         \vdots & \vdots & \ddots & \vdots \\
         0 & 0 & \cdots & 1 \\
         -\alpha_{0} & -\alpha_{1} & \cdots & -\alpha_{l-1} \\
        \end{smallmatrix}
    \end{bmatrix}$ and
$\Gamma=\begin{bmatrix}
        \begin{smallmatrix}
         1 \\
          0 \\
           \vdots \\
            0 \\
        \end{smallmatrix}
    \end{bmatrix}^{T}$,
where  $\alpha_{0}, \alpha_{1},  \cdots,  \alpha_{l-1}$ are some real numbers such that $m(\mathbf{s})=\mathbf{s}^{l}+\alpha_{l-1}\mathbf{s}^{l-1}+\cdots+\alpha_{1}\mathbf{s}+\alpha_{0}$ is the minimal polynomial of $S$.
Let $T$  be the unique solution of the following
Sylvester equation \cite{nik98}:
\begin{equation}\label{Sylvester}
\begin{split}
  T\Phi-MT=Q\Gamma,
\end{split}
\end{equation}
and,  for $i=1,\cdots,N$, let $\Pi_{i}(w)$ be the unique solution
of the following Sylvester equation
\begin{equation}\label{Pi}
\begin{split}
\Pi_{i}(w)S\!=\!A_{1i}(w)\Pi_{i}(w)\!+\!A_{2i}(w)F(w)\!+\!E_{0i}(w).
\end{split}
\end{equation}
Finally, let $U_{i}(w)=-b_{i}^{-1}(w)\big(A_{3i}(w)\Pi_{i}(w)+E_{ri}(w)+\sum_{s=1}^{r}c_{si}(w)F(w)S^{s-1}-F(w)S^{r}\big)$
and  $\Upsilon_{i}(w)=\mbox{col}(U_{i}(w),U_{i}(w)S$, $\cdots,U_{i}(w)S^{l-1})$.

Then, as shown in \cite{Su6}, performing  on the internal model \eqref{dotetai} and the plant \eqref{system1} the  following transformation
\begin{equation}\label{transformation1}
\begin{split}
  &\bar{z}_{i}=z_{i}-\Pi_{i}(w)v\\
  & \bar{\xi}_{si}=\xi_{si}-F(w)S^{s-1}v,~ s=1,\cdots,r\\
 & \bar{\eta}_{i}=\eta_{i}-T\Upsilon_{i}(w)v\\
  & \bar{u}_{i}=u_{i}-\Gamma T^{-1}\eta_{i},~ i=1,\cdots,N
\end{split}
\end{equation}
leads to the following so-called augmented system
\begin{equation}\label{system4}
\begin{split}
  \dot{\bar{z}}_{i}=&A_{1i}(w)\bar{z}_{i}+A_{2i}(w)\bar{\xi}_{1i}\\
  \dot{\bar{\xi}}_{si}=&\bar{\xi}_{(s+1)i},\ s=1,\cdots,r-1\\
  \dot{\bar{\xi}}_{ri}=&A_{3i}(w)\bar{z}_{i}+\sum_{s=1}^{r}c_{si}(w)\bar{\xi}_{si}+b_{i}(w)\bar{u}_{i}\\
  &+b_{i}(w)\Gamma T^{-1}\bar{\eta}_{i}\\
  \dot{\bar{\eta}}_{i}=&(M+Q\Gamma T^{-1})\bar{\eta}_{i}+Q\bar{u}_{i},\\
  e_{i}=&\bar{\xi}_{1i}.
\end{split}
\end{equation}


Then  the main result in \cite{Su6} can be summarized as follows.

\begin{Lemma}\label{Lemma1}
Consider the following continuous-time distributed output feedback control law 
\begin{equation}\label{bui2}
\begin{split}
\bar{u}_{i}(t)&=-K\Omega\zeta_{i}(t)\\
\dot{\zeta}_{i}(t)&=A_{o}(h)\zeta_{i}(t)+B_{o}(h)e_{vi}(t),\ \ i=1,\cdots,N
\end{split}
\end{equation}
where $\zeta_{i}=\mbox{col}(\zeta_{1i},\zeta_{2i},\cdots,\zeta_{ri})\in\mathbb{R}^{r}$, $\Omega=[\gamma_{0},\gamma_{1},\cdots$, $\gamma_{r-2},1]$,
$A_{o}(h)\!=\!\!\!\begin{bmatrix}
        \begin{smallmatrix}
              -h\delta_{r-1} & 1 & 0 & \cdots & 0 \\
              -h^{2}\delta_{r-2} & 0 & 1 & \cdots & 0 \\
              \vdots & \vdots & \vdots & \ddots & \vdots \\
              -h^{r-1}\delta_{1} & 0 & 0 & \cdots & 1 \\
              -h^{r}\delta_{0} & 0 & 0 & 0 & 0 \\
        \end{smallmatrix}
    \end{bmatrix}$,
$B_{o}(h)\!=\!\!\!\begin{bmatrix}
        \begin{smallmatrix}
                                   h\delta_{r-1} \\
                                   h^{2}\delta_{r-2} \\
                                   \vdots \\
                                   h^{r-1}\delta_{1} \\
                                   h^{r}\delta_{0} \\
        \end{smallmatrix}
    \end{bmatrix}$,
$h$ and $K$ are some positive real numbers,
and $\gamma_{i}$, $i=0,1,\cdots,r-2$, and $\delta_{j}$, $j=0,1,\cdots,r-1$, are some real numbers such that the polynomials $f(\mathbf{s})=\mathbf{s}^{r-1}+\gamma_{r-2}\mathbf{s}^{r-2}+\cdots+\gamma_{1}\mathbf{s}+\gamma_{0}$ and  $g(\mathbf{s})=\mathbf{s}^{r}+\delta_{r-1}\mathbf{s}^{r-1}+\cdots+\delta_{1}\mathbf{s}+\delta_{0}$ are both stable. Under Assumptions \ref{Ass1}-\ref{Ass3},
there exist some positive real numbers $K^{*}$ and $h^{*}$ such that, for any $K>K^{*}$ and $h>h^{*}$, 
the  control law (\ref{bui2}) solves the cooperative robust stabilization problem of the system \eqref{system4}.
\end{Lemma}

\begin{Remark} As noted in \cite{Su6}, 
  if  the continuous-time distributed output feedback control law \eqref{bui2} solves the cooperative robust stabilization problem of the augmented system (\ref{system4}), 
  then the following continuous-time distributed output feedback control law
\begin{equation}\label{ui3}
\begin{split}
u_{i}(t)&=\Gamma T^{-1}\eta_{i}(t)-K\Omega\zeta_{i}(t)\\
\dot{\eta}_{i}(t)&=M\eta_{i}(t)+Q(\Gamma T^{-1}\eta_{i}(t)-K\Omega\zeta_{i}(t))\\
\dot{\zeta}_{i}(t)&=A_{o}(h)\zeta_{i}(t)+B_{o}(h)e_{vi}(t),\ \ i=1,\cdots,N
\end{split}
\end{equation}
solves the cooperative robust output regulation problem of the system (\ref{system1}).
\end{Remark}


\section{Main Result}\label{MR}
In this section, we will consider the cooperative robust practical output regulation problem for the system \eqref{system1} by a distributed output-based event-triggered control law.

In order to obtain a digital control law, for any $t\in[t_{k}^{i},t_{k+1}^{i})$ with $k\in\mathbb{Z}^{+}$ and $i=1,\cdots,N$,
we perform a  new transformation instead of \eqref{transformation1} as follows
\begin{equation}\label{transformation2}
\begin{split}
 & \bar{z}_{i}=z_{i}-\Pi_{i}(w)v \\
 & \bar{\xi}_{si}=\xi_{si}-F(w)S^{s-1}v,~ s=1,\cdots,r\\
 & \bar{\eta}_{i}=\eta_{i}-T\Upsilon_{i}(w)v\\
 & \hat{u}_{i}=u_{i}-\Gamma T^{-1}\eta_{i}(t_{k}^{i}).\\
\end{split}
\end{equation}
Then we get a new augmented system
\begin{equation}\label{system4a}
\begin{split}
  \dot{\bar{z}}_{i}=&A_{1i}(w)\bar{z}_{i}+A_{2i}(w)\bar{\xi}_{1i}\\
  \dot{\bar{\xi}}_{si}=&\bar{\xi}_{(s+1)i},\ s=1,\cdots,r-1\\
  \dot{\bar{\xi}}_{ri}=&A_{3i}(w)\bar{z}_{i}+\sum_{s=1}^{r}c_{si}(w)\bar{\xi}_{si}+b_{i}(w)\hat{u}_{i}\\
  &+b_{i}(w)\Gamma T^{-1}(\bar{\eta}_{i}+\tilde{\eta}_{i})\\
  \dot{\bar{\eta}}_{i}=&(M+Q\Gamma T^{-1})\bar{\eta}_{i}+Q\hat{u}_{i}+Q\Gamma T^{-1}\tilde{\eta}_{i}\\
  e_{i}=&\bar{\xi}_{1i}
\end{split}
\end{equation}
for any $t\in[t_{k}^{i},t_{k+1}^{i})$ with $k\in\mathbb{Z}^{+}$ and $i=1,\cdots,N$.

Following the framework in \cite{Su6}, the first step here is also to stabilize the augmented system \eqref{system4a}. For this purpose, according to the continuous-time distributed output feedback control law \eqref{bui2}, we consider the following piecewise constant distributed output-based event-triggered  control law
\begin{equation}\label{bui3}
\begin{split}
\hat{u}_{i}(t)&=-K\Omega\zeta_{i}(t_{k}^{i})\\
\dot{\zeta}_{i}(t)&=A_{o}(h)\zeta_{i}(t)+B_{o}(h)e_{vi}(t_{k}^{i})
\end{split}
\end{equation}
for any $t\in[t_{k}^{i},t_{k+1}^{i})$ with $k\in\mathbb{Z}^{+}$ and $i=1,\cdots,N$.
Thus we obtain the closed-loop system composed of the augmented system \eqref{system4a} and the control law \eqref{bui3} as follows
\begin{equation}\label{xc1}
\begin{split}
  \dot{\bar{z}}_{i}=&A_{1i}(w)\bar{z}_{i}+\bar{A}_{2i}(w)\bar{\xi}_{i}\\
  \dot{\bar{\xi}}_{i}\!=&A_{4i}(w)\bar{z}_{i}\!+\!A_{5i}(w)\bar{\xi}_{i}\!+\!B_{0i}(w)\Gamma T^{-1}(\bar{\eta}_{i}+\tilde{\eta}_{i})\\
  &-B_{0i}(w)K\Omega(\zeta_{i}+\tilde{\zeta}_{i})\\
  \dot{\bar{\eta}}_{i}\!=&(M+Q\Gamma T^{-1})\bar{\eta}_{i}-QK\Omega(\zeta_{i}\!+\!\tilde{\zeta}_{i})\!+\!Q\Gamma T^{-1}\tilde{\eta}_{i}\\
  \dot{\zeta}_{i}=&A_{o}(h)\zeta_{i}+B_{o}(h)(e_{vi}+\tilde{e}_{vi}),~ i=1,\cdots,N
\end{split}
\end{equation}
where $\bar{\xi}_{i}=\mbox{col}(\bar{\xi}_{1i},\cdots,\bar{\xi}_{ri})$,
$\bar{A}_{2i}(w)=A_{2i}(w)D,~ D=[1\ 0\ \cdots\ 0]_{1\times r}$,
$A_{4i}(w)=\begin{bmatrix}
        \begin{smallmatrix}
          0_{(r-1)\times(n_{i}-r)} \\
          A_{3i}(w) \\
        \end{smallmatrix}
    \end{bmatrix}$,
$B_{0i}(w)=[0\ \cdots\ 0\ b_{i}(w)]^{T}_{1\times r}$,
$A_{5i}(w)=\begin{bmatrix}
        \begin{smallmatrix}
                  0 & 1 & \cdots & 0 \\
                  \vdots & \vdots & \ddots & \vdots \\
                  0 & 0 & \cdots & 1 \\
                  c_{1i}(w) & c_{2i}(w) & \cdots & c_{ri}(w) \\
        \end{smallmatrix}
    \end{bmatrix}$,
Let $\bar{z}=\mbox{col}(\bar{z}_{1},\cdots,\bar{z}_{N})$, $\bar{\xi}=\mbox{col}(\bar{\xi}_{1},\cdots,\bar{\xi}_{N})$, $\bar{\eta}=\mbox{col}(\bar{\eta}_{1}$, $\cdots,\bar{\eta}_{N})$, $\tilde{\eta}=\mbox{col}(\tilde{\eta}_{1},\cdots,\tilde{\eta}_{N})$, $\zeta=\mbox{col}(\zeta_{1},\cdots,\zeta_{N})$ and $\tilde{\zeta}=\mbox{col}(\tilde{\zeta}_{1},\cdots,\tilde{\zeta}_{N})$. Then the system \eqref{xc1} can be put into the following form
\begin{equation}\label{xc2}
\begin{split}
  \dot{\bar{z}}=&\bar{A}_{1}(w)\bar{z}+\bar{A}_{2}(w)\bar{\xi}\\
  \dot{\bar{\xi}}=&\bar{A}_{4}(w)\bar{z}\!+\!\bar{A}_{5}(w)\bar{\xi}\!+\!\bar{A}_{6}(w)(\bar{\eta}\!+\!\tilde{\eta})\\
  &+\bar{A}_{7}(w)(\zeta+\tilde{\zeta})\\
  \dot{\bar{\eta}}=&\bar{A}_{8}\bar{\eta}+\bar{A}_{9}(\zeta+\tilde{\zeta})+\bar{A}_{10}\tilde{\eta}\\
  \dot{\zeta}=&\bar{A}_{11}\bar{\xi}+\bar{A}_{12}\zeta+\tilde{e}_{v}(t)\otimes B_{o}(h)
\end{split}
\end{equation}
where
$\bar{A}_{2}(w)=\mbox{blockdiag}\{\bar{A}_{21}(w),\cdots,\bar{A}_{2N}(w)\}$,
$\bar{A}_{k}(w)=\mbox{blockdiag}\{A_{k1}(w),\cdots,A_{kN}(w)\},~k\!=\!1,4,5$,
$\bar{A}_{6}(w)= \mbox{blockdiag}\{B_{01}(w)\Gamma T^{-1},\cdots,B_{0N}(w)\Gamma T^{-1}\}$,
$\bar{A}_{7}(w)=\mbox{blockdiag}\{-B_{01}(w)K\Omega,\cdots,-B_{0N}(w)K\Omega\}$,
$\bar{A}_{8}= I_{N}\otimes(M+Q\Gamma T^{-1}),~\bar{A}_{9}=-I_{N}\otimes (QK\Omega)$,
$\bar{A}_{10}= I_{N}\otimes(Q\Gamma T^{-1}),~\bar{A}_{12}=I_{N}\otimes A_{o}(h)$,
$\bar{A}_{11}=(H\otimes I_{r})\bar{D},~\bar{D}=I_{N}\otimes[B_{o}(h)\ 0_{r\times(r-1)}]$.
Let $x_{c}=\mbox{col}(z,\xi,\eta,\zeta)$ and $\bar{x}_{c}=\mbox{col}(\bar{z},\bar{\xi},\bar{\eta},\zeta)$. Then the system \eqref{xc2} can be further put into the following compact form
\begin{equation}\label{xc3}
\begin{split}
\dot{\bar{x}}_{c}=A_{c}(w)\bar{x}_{c}+f_{c}(\tilde{\eta},\tilde{\zeta},\tilde{e}_{v})
\end{split}
\end{equation}
with
$A_{c}(w)=\begin{bmatrix}
        \begin{smallmatrix}
                        \bar{A}_{1}(w) & \bar{A}_{2}(w) & 0 & 0 \\
                        \bar{A}_{4}(w) & \bar{A}_{5}(w) & \bar{A}_{6}(w) & \bar{A}_{7}(w) \\
                        0 & 0 & \bar{A}_{8}& \bar{A}_{9} \\
                        0 & \bar{A}_{11} & 0 & \bar{A}_{12} \\
        \end{smallmatrix}
    \end{bmatrix}$,
    $f_{c}(\tilde{\eta},\tilde{\zeta},\tilde{e}_{v})=\begin{bmatrix}
        \begin{smallmatrix}
                              0 \\
                              \bar{A}_{6}(w)\tilde{\eta}+\bar{A}_{7}(w)\tilde{\zeta} \\
                              \bar{A}_{9}\tilde{\zeta}+\bar{A}_{10}\tilde{\eta} \\
                              \tilde{e}_{v}(t)\otimes B_{o}(h) \\
        \end{smallmatrix}
    \end{bmatrix}$.

\begin{Proposition}\label{Proposition1}
Under Assumptions \ref{Ass1}-\ref{Ass3}, for any $\epsilon>0$, and any known compact sets $\mathbb{V}$ and $\mathbb{W}$, if the solution $\bar{x}_{c}(t)$ of the closed-loop system \eqref{xc3} exists and is bounded for all $t\geq0$, and satisfies
 \begin{equation}\label{bxc1}
\begin{split}
\lim_{t \to +\infty} \sup \|\bar{x}_{c}(t)\|\leq\epsilon,
\end{split}
\end{equation}
then Problem \ref{LCORPS} for the system \eqref{system1} is solvable by the following distributed event-triggered  output feedback control law
 \begin{equation}\label{ui4}
\begin{split}
u_{i}(t)&=\Gamma T^{-1}\eta_{i}(t_{k}^{i})-K\Omega\zeta_{i}(t_{k}^{i})\\
\dot{\eta}_{i}(t)&=M\eta_{i}(t)+Q(\Gamma T^{-1}\eta_{i}(t_{k}^{i})-K\Omega\zeta_{i}(t_{k}^{i}))\\
\dot{\zeta}_{i}(t)&=A_{o}(h)\zeta_{i}(t)+B_{o}(h)e_{vi}(t_{k}^{i})
\end{split}
\end{equation}
for $t\in[t_{k}^{i},t_{k+1}^{i})$ with $k\in\mathbb{Z}^{+}$ and $i=1,\cdots,N$.
\end{Proposition}
\begin{Proof}
First, it is easy to see that \eqref{bxc1} implies that
\begin{equation}
\begin{split}
\lim_{t \to +\infty} \sup \|e(t)\|\leq\lim_{t \to +\infty} \sup \|\bar{x}_{c}(t)\|\leq\epsilon.
\end{split}
\end{equation}
Second, according to the transformation \eqref{transformation2}, we have
$ z_{i}= \bar{z}_{i}+\Pi_{i}(w)v$, $\xi_{si}=\bar{\xi}_{si}+F(w)S^{s-1}v$, $s=1,\cdots,r$,   $\eta_{i}=\bar{\eta}_{i}+T\Upsilon_{i}(w)v$ for $i=1,\cdots,N$.
 Since $v(t)$ and $w$ are bounded for all $t\geq0$,   $\Pi_{i}(w)v$, $F(w)S^{s-1}v$ and  $T\Upsilon_{i}(w)v$ for  $s=1,\cdots,r$ and $i=1,\cdots,N$ are also bounded for all $t\geq0$. Thus,  the fact that $\bar{x}_{c}(t)$ exists and is bounded for all $t\geq0$ implies that $x_{c}(t)$ also exists and is bounded for all $t\geq0$.

 That is to say, the two properties of Problem \ref{LCORPS} are both satisfied.
 \end{Proof}

\begin{Remark}\label{RemarkControl1}
It is easy to see that the control law \eqref{ui4} takes the same form as \eqref{ui1}. As noted in Remark \ref{RemarkProb1a}, we can further discretize the internal model and the dynamic compensator and get a digital control law of the form \eqref{ui2}.
\end{Remark}
We say the distributed control law \eqref{bui3} solves the cooperative robust practical stabilization problem for \eqref{system4a} if, for any given
$\epsilon >0$,
 $\bar{x}_{c}(t)$ exists and is bounded for all $t\geq0$, and satisfies \eqref{bxc1}. For simplicity, we will use the abbreviation CRPSP to denote cooperative robust practical stabilization problem.

According to Lemma \ref{Lemma1}, 
 there always exist some positive numbers $K^{*}$ and $h^{*}$ such that for any $K>K^{*}$ and $h>h^{*}$, the closed-loop system matrix $A_{c} (w) $ is  Hurwitz for all $w\in\mathbb{W}$. Then the following Lyapunov equation
 \begin{equation}\label{LMI2}
\begin{split}
A_{c}(w)^{T}P(w)+P(w)A_{c}(w)=-2I
  \end{split}
\end{equation}
always has a solution $P(w)$ which is  symmetric and positive definite.

 Next, we introduce some notation as follows. Let $\lambda_{M}=\sup_{w\in\mathbb{W}}\{\lambda_{\max}(P(w))\}$ and $\lambda_{m}=\inf_{w\in\mathbb{W}}\{\lambda_{\min}(P(w))\}$.  Let $\lambda_{1}=\sup_{w\in\mathbb{W}}\{\|P(w)\|^{2}\}$, $\lambda_{2}=\max_{i=1,...,\rho}\{\|H_{i}\|^{2}\}$,  where $\rho$ is the total number of connected directed graphs associated with $N+1$ agents, $\lambda_{3}=\|B_{o}(h)\|^{2}$, $\lambda_{4}=\sup_{w\in\mathbb{W}}\{\|\bar{A}_{6}(w)\|^{2}+\|\bar{A}_{10}\|^{2}\}$, $\lambda_{5}=\sup_{w\in\mathbb{W}}\{\|\bar{A}_{7}(w)\|^{2}+\|\bar{A}_{9}\|^{2}\}$, $\bar{\lambda}_{2}=\max\{1,\lambda_{2}\}$ and $\bar{\lambda}_{M}=\max\{\lambda_{3},\lambda_{4},\lambda_{5}\}$. Let $v_{ci}=\mbox{col}(\eta_{i},\zeta_{i},e_{vi})$,
 $\tilde{v}_{ci}=\mbox{col}(\tilde{\eta}_{i},\tilde{\zeta}_{i},\tilde{e}_{vi})$ for $i=1,\cdots,N$, $v_{c}=\mbox{col}(v_{c1},\cdots,v_{cN})$ and $\tilde{v}_{c}=\mbox{col}(\tilde{v}_{c1},\cdots,\tilde{v}_{cN})$. Then we introduce the following distributed output-based event-triggered mechanism
\begin{equation}\label{trigger2}
\begin{split}
t_{k+1}^{i}\!=\!\inf\{t>t_{k}^{i}~|~ \|\tilde{v}_{ci}(t)\|^{2}\geq& \bar{\sigma}(\|\zeta_{i}(t)\|^{2}+e_{vi}^{2}(t))\\
&+\!\beta \mathbf{e}^{-\alpha t}+\delta\}\\
\end{split}
\end{equation}
where $\bar{\sigma}=\frac{\sigma}{\lambda_{1}\bar{\lambda}_{M}\bar{\lambda}_{2}}$, $\sigma\in(0,1)$, $\beta>0$, $\alpha\in(0,\lambda)$ with $\lambda=\frac{1-\sigma}{\lambda_{M}}$ and $\delta>0$ are some constants. Under the event-triggered mechanism \eqref{trigger2}, we easily obtain that, for any $t\in[t_{k}^{i},t_{k+1}^{i})$  with $k\in\mathbb{Z}^{+}$ and $i=1,\cdots,N$,
\begin{equation}\label{tildeevi1}
\begin{split}
\|\tilde{v}_{ci}(t)\|^{2}\!\leq\! \frac{\sigma}{\lambda_{1}\bar{\lambda}_{M}\bar{\lambda}_{2}}(\|\zeta_{i}(t)\|^{2}\!+\!e_{vi}^{2}(t))\!+\!\beta \mathbf{e}^{-\alpha t}\!+\!\delta.
\end{split}
\end{equation}
\begin{Remark}
%
%
%
It should be noted that here the exponential convergence term $\beta\mathbf{e}^{-\alpha t}$ governs the triggering threshold and influences the event-triggered numbers during the initial stage, while the constant $\delta$ governs the triggering threshold and influences the event-triggered numbers during the steady stage.
\end{Remark}

Let $[0,T_{M})$ with $0<T_{M}\leq\infty$ denote the right maximally defined solution interval of the closed-loop system \eqref{xc3} under the event-triggered mechanism \eqref{trigger2}.
Then we have the following lemma.
\begin{Lemma}\label{Lemma2}
Under Assumptions \ref{Ass1}-\ref{Ass3}, for any compact sets $\mathbb{V} \subset \mathbb{R}^{n_{v}}$ and $\mathbb{W} \subset \mathbb{R}^{n_{w}} $ with $0\in\mathbb{V}$ and $0\in\mathbb{W}$,
there exist some positive real numbers $K^{*}$ and $h^{*}$ such that, for any $K>K^{*}$ and $h>h^{*}$,
\begin{equation}\label{bxct1}
\begin{split}
\|\bar{x}_{c}(t)\|^{2}&\leq a\mathbf{e}^{-\alpha t} +b,~\forall t\in[0,T_{M})\\
\end{split}
\end{equation}
where $a$ and $b$  are some positive constants to be specified later. As a result, $\bar{x}_{c}(t)$ is bounded for all $t\in[0,T_{M})$.
\end{Lemma}
\begin{Proof} 
 First, it is ready to verify that
\begin{equation}\label{fc1}
\begin{split}
\|f_{c}(\tilde{\eta},\tilde{\zeta},\tilde{e}_{v})\|^{2}\!\!=&\|\bar{A}_{9}\tilde{\zeta}\!+\!\bar{A}_{10}\tilde{\eta}\|^{2}\!+\!\|\tilde{e}_{v}(t)\otimes B_{o}(h)\|^{2} \\
&+\|\bar{A}_{6}(w)\tilde{\eta}+\bar{A}_{7}(w)\tilde{\zeta}\|^{2}\\
\leq& (\|\bar{A}_{6}(w)\|^{2}+\|\bar{A}_{10}\|^{2})\|\tilde{\eta}\|^{2}\\
&+(\|\bar{A}_{7}(w)\|^{2}+\|\bar{A}_{9}\|^{2})\|\tilde{\zeta}\|^{2}\\
&+\|B_{o}(h)\|^{2}\|\tilde{e}_{v}\|^{2}\\
\leq&\lambda_{3}\|\tilde{e}_{v}\|^{2}+\lambda_{4}\|\tilde{\eta}\|^{2}+\lambda_{5}\|\tilde{\zeta}\|^{2}\\
\leq&\bar{\lambda}_{M}\|\tilde{v}_{c}\|^{2}.
\end{split}
\end{equation}
By \eqref{tildeevi1} and the definition of  $\bar{\lambda}_{2}$, 
we have
\begin{equation}\label{tildevc1}
\begin{split}
&\|\tilde{v}_{c}(t)\|^{2}=\sum_{i=1}^{N}\|\tilde{v}_{ci}(t)\|^{2}\\
\leq&\sum_{i=1}^{N}\frac{\sigma}{\lambda_{1}\bar{\lambda}_{M}\bar{\lambda}_{2}}(\|\zeta_{i}(t)\|^{2}\!+\!e_{vi}^{2}(t))+N\beta \mathbf{e}^{-\alpha t}+N\delta\\
=&\frac{\sigma}{\lambda_{1}\bar{\lambda}_{M}\bar{\lambda}_{2}}(\|\zeta(t)\|^{2}+\|e_{v}(t)\|^{2})+N\beta \mathbf{e}^{-\alpha t}+N\delta\\
\leq&\frac{\sigma}{\lambda_{1}\bar{\lambda}_{M}\bar{\lambda}_{2}}(\|\zeta(t)\|^{2}\!+\!\|H\|^{2}\|e(t)\|^{2})\!+\!\!N\beta \mathbf{e}^{-\alpha t}\!+\!\!N\delta\\
\leq&\frac{\sigma}{\lambda_{1}\bar{\lambda}_{M}}(\|\zeta(t)\|^{2}+\|e(t)\|^{2})+N\beta \mathbf{e}^{-\alpha t}+N\delta\\
\leq& \frac{\sigma}{\lambda_{1}\bar{\lambda}_{M}}\|\bar{x}_{c}(t)\|^{2}+N\beta \mathbf{e}^{-\alpha t}+N\delta\\
\end{split}
\end{equation}
which implies that
\begin{equation}\label{tildevc2}
\begin{split}
\lambda_{1}\bar{\lambda}_{M}\|\tilde{v}_{c}(t)\|^{2}\leq& \sigma\|\bar{x}_{c}(t)\|^{2}+\lambda_{1}\bar{\lambda}_{M}N\beta \mathbf{e}^{-\alpha t}\\
&+\lambda_{1}\bar{\lambda}_{M}N\delta.\\
\end{split}
\end{equation}
Also, according to \eqref{fc1} and the definition of  $\lambda_{1}$, we have
\begin{equation}\label{PB1}
\begin{split}
\|P(w)\|^{2}\|f_{c}(\tilde{\eta},\tilde{\zeta},\tilde{e}_{v})\|^{2}\leq \lambda_{1}\bar{\lambda}_{M}\|\tilde{v}_{c}\|^{2}.\\
\end{split}
\end{equation}
Combining \eqref{tildevc2} and \eqref{PB1} gives
\begin{equation}\label{tildevc3}
\begin{split}
\|P(w)\|^{2}\|f_{c}(\tilde{\eta},\tilde{\zeta},\tilde{e}_{v}\!)\|^{2}\!\!\leq\! \sigma\|\bar{x}_{c}(t)\|^{2}\!\!+\!\lambda_{0}\beta \mathbf{e}^{-\alpha t}\!\!+\!\lambda_{0}\delta\\
\end{split}
\end{equation}
where $\lambda_{0}=\lambda_{1}\bar{\lambda}_{M}N$.

Let $V(\bar{x}_{c}(t))=\bar{x}_{c}(t)^{T}P(w)\bar{x}_{c}(t)$. Then
\begin{equation}\label{Vxc1}
\begin{split}
\lambda_{m}\|\bar{x}_{c}(t)\|^{2}\leq V(\bar{x}_{c}(t))\leq \lambda_{M}\|\bar{x}_{c}(t)\|^{2}.\\
\end{split}
\end{equation}
According to  \eqref{tildevc3} and \eqref{Vxc1}, we have 
\begin{equation}\label{dotVxc1}
\begin{split}
&\frac{\partial V(\bar{x}_{c}) }{\partial \bar{x}_{c}}(A_{c}(w)\bar{x}_{c}+f_{c}(\tilde{\eta},\tilde{\zeta},\tilde{e}_{v}))\\
=&\bar{x}_{c}^{T}(t)(A_{c}^{T}(w)P(w)\!+\!P(w)A_{c}(w))\bar{x}_{c}(t)\\
&+2\bar{x}_{c}^{T}(t)P(w)f_{c}(\tilde{\eta},\tilde{\zeta},\tilde{e}_{v})\\
\leq& -2\|\bar{x}_{c}(t)\|^{2}+\|\bar{x}_{c}(t)\|^{2}\\
&+\|P(w)\|^{2}\|f_{c}(\tilde{\eta},\tilde{\zeta},\tilde{e}_{v})\|^{2}\\
\leq& -\|\bar{x}_{c}(t)\|^{2}+\sigma\|\bar{x}_{c}(t)\|^{2}+\lambda_{0}\beta \mathbf{e}^{-\alpha t}\!+\!\lambda_{0}\delta\\
\leq &-\frac{1-\sigma}{\lambda_{M}}V(\bar{x}_{c}(t))\!+\!\lambda_{0}\beta \mathbf{e}^{-\alpha t}\!+\!\lambda_{0}\delta\\
=&-\lambda V(\bar{x}_{c}(t))\!+\!\lambda_{0}\beta \mathbf{e}^{-\alpha t}\!+\!\lambda_{0}\delta.\\
\end{split}
\end{equation}
By \eqref{dotVxc1}, we can obtain
\begin{equation}\label{Vxc2}
\begin{split}
V(\bar{x}_{c}(t))
\leq&\mathbf{e}^{-\lambda t}V(\bar{x}_{c}(0))+\frac{\lambda_{0}\beta}{\lambda-\alpha}\mathbf{e}^{-\alpha t}+\frac{\lambda_{0}\delta}{\lambda}\\
\leq& (V(\bar{x}_{c}(0))+\frac{\lambda_{0}\beta}{\lambda-\alpha})\mathbf{e}^{-\alpha t} +\frac{\lambda_{0}\delta}{\lambda}.\\
\end{split}
\end{equation}
By \eqref{Vxc1} and \eqref{Vxc2}, we have
\begin{equation}\label{xct1}
\begin{split}
\|\bar{x}_{c}(t)\|^{2}&\!\leq\frac{V(\bar{x}_{c}(t))}{\lambda_{m}}\\
&\!\leq(\frac{V(\bar{x}_{c}(0))}{\lambda_{m}}\!+\!\frac{\lambda_{0}\beta}{\lambda_{m}(\lambda\!-\!\alpha)})\mathbf{e}^{\!-\alpha t}\!+\!\frac{\lambda_{0}\delta}{\lambda\lambda_{m}} \\
&=a\mathbf{e}^{-\alpha t}+b
\end{split}
\end{equation}
with $a=\frac{V(\bar{x}_{c}(0))}{\lambda_{m}}\!+\!\frac{\lambda_{0}\beta}{\lambda_{m}(\lambda\!-\!\alpha)}$ and $b=\frac{\lambda_{0}\delta}{\lambda\lambda_{m}}$.

From \eqref{xct1}, we have
\begin{equation}\label{bxct2}
\begin{split}
\|\bar{x}_{c}(t)\|\leq\sqrt{a\mathbf{e}^{-\alpha t}+b}\leq\sqrt{a+b},~\forall t\in[0,T_{M}).
\end{split}
\end{equation}
That is to say, $\bar{x}_{c}(t)$ is bounded for all $t\in[0,T_{M})$.
Thus the proof is completed.
\end{Proof}

Combining Proposition \ref{Proposition1} and  Lemma \ref{Lemma2} gives
the following main result.
\begin{Theorem}\label{Theorem1}
Under Assumptions \ref{Ass1}-\ref{Ass3}, for any compact sets $\mathbb{V} \subset \mathbb{R}^{n_{v}}$ and $\mathbb{W} \subset \mathbb{R}^{n_{w}} $ with $0\in\mathbb{V}$ and $0\in\mathbb{W}$, there exist some positive real numbers $K^{*}$ and $h^{*}$ such that, for any $K>K^{*}$ and $h>h^{*}$,  the distributed output feedback control law \eqref{ui4}
 together with the distributed output-based event-triggered mechanism \eqref{trigger2} with $\delta=\frac{\lambda\lambda_{m}}{\lambda_{0}}\epsilon^{2}$ solves the cooperative robust practical output regulation problem for the system \eqref{system1}.
\end{Theorem}
\begin{Proof}
By Proposition \ref{Proposition1} and  Lemma \ref{Lemma2}, it suffices to show that $T_M = \infty$. For this purpose,
let us first consider the case where, for all $i = 1, \cdots, N$, the number of the triggering times is finite. Then, there exists some finite time $T_0$ such that the closed-loop system is a linear-time-invariant continuous-time system for $t \geq T_0$. Thus $T_M$ must be equal to infinity. Next, we show that,  for any $i = 1, \cdots, N$, if the sequence $\{t_{k}^{i} \}$ has infinite many members, then
$\lim_{k \rightarrow \infty} t_{k}^{i} = \infty$. Thus, $T_M = \infty$.

In fact, by \eqref{trigger2}, we have
\begin{equation}\label{tildevci1}
\begin{split}
\|\tilde{v}_{ci}((t_{k+1}^{i})^{-})\|
&\geq \sqrt{\delta}.\\
\end{split}
\end{equation}
According to \eqref{tildevc1} and \eqref{bxct2}, for any $t\in[0,T_{M})$, we have
\begin{equation}\label{tildevc4}
\begin{split}
\|\tilde{v}_{c}(t)\|^{2}\leq& \frac{\sigma}{\lambda_{1}\bar{\lambda}_{M}}\|\bar{x}_{c}(t)\|^{2}+N\beta \mathbf{e}^{-\alpha t}+N\delta\\
\leq&\frac{\sigma}{\lambda_{1}\bar{\lambda}_{M}}(a+b)+N\beta+N\delta.
\end{split}
\end{equation}
Based on \eqref{xc3} and \eqref{fc1}, for any $t\in[0,T_{M})$, we have 
\begin{equation*}
\begin{split}
&\frac{d}{d t} \|\tilde{v}_{ci}(t)\|=\frac{d}{d t} (\tilde{v}_{ci}^{T}(t)\tilde{v}_{ci}(t))^{\frac{1}{2}}\\
=&(\tilde{v}_{ci}^{T}(t)\tilde{v}_{ci}(t))^{-\frac{1}{2}}\tilde{v}_{ci}^{T}(t)\dot{\tilde{v}}_{ci}(t)\\
\leq&\|\dot{\tilde{v}}_{ci}(t)\|=\|-\dot{v}_{ci}(t)\|=\|\mbox{col}(\dot{\eta}_{i},\dot{\zeta}_{i},\dot{e}_{vi})\|\\
=&\|\mbox{col}(\dot{\bar{\eta}}_{i}+T\Upsilon_{i}(w)\dot{v},\dot{\zeta}_{i},\sum_{j=1}^{N}h_{ij}\dot{e}_{j}(t))\|\\
\end{split}
\end{equation*}
\begin{equation}\label{dottildevci1}
\begin{split}
\leq&\|T\Upsilon_{i}(w)Sv\|+N\|\dot{\bar{x}}_{c}(t)\|\\
=&\|T\Upsilon_{i}(w)Sv\|+N\|A_{c}(w)\bar{x}_{c}+f_{c}(\tilde{\eta},\tilde{\zeta},\tilde{e}_{v})\|\\
\leq&\|T\Upsilon_{i}(w)Sv\|\!\!+\!\!N\|A_{c}(w)\|\|\bar{x}_{c}(t)\|\\
&+N\|f_{c}(\tilde{\eta},\tilde{\zeta},\tilde{e}_{v})\|\\
\leq&\|T\Upsilon_{i}(w)Sv\|+N\|A_{c}(w)\|\|\bar{x}_{c}(t)\|\\
&+N\sqrt{\bar{\lambda}_{M}}\|\tilde{v}_{c}(t)\|.\\
\end{split}
\end{equation}
Note that, for any $v(t)\in\mathbb{V}$ and any $w\in\mathbb{W}$ with $\mathbb{V}$ and $\mathbb{W}$ being some compact sets,  $\|T\Upsilon_{i}(w)Sv\|$ and $A_{c}(w)$ are bounded. Also, according to \eqref{bxct2} and \eqref{tildevc4},  $\bar{x}_{c}(t)$ and $\tilde{v}_{c}(t)$ are bounded for all $t\in[0,T_{M})$. Thus, for $i=1,\cdots,N$,  all $v(t)\in\mathbb{V}$, all $w\in\mathbb{W}$, and all $\bar{x}_{c}(0)\in\mathbb{\bar{X}}$ with $\mathbb{V}$, $\mathbb{W}$ and $\mathbb{\bar{X}}$ being some compact sets, there always exists a positive constant $c$ depending on $\delta$ such that
\begin{equation}\label{dottildevci2}
\begin{split}
\bigg|\frac{d}{d t} \|\tilde{v}_{ci}(t)\|\bigg|\leq c,~\forall t\in[0,T_{M}).
\end{split}
\end{equation}
%
Combining \eqref{tildevci1} and \eqref{dottildevci2}, we have
\begin{equation}\label{tki3}
\begin{split}
t_{k+1}^{i}-t_{k}^{i}&\geq \frac{\tilde{v}_{ci}((t_{k+1}^{i})^{-})-\tilde{v}_{ci}(t_{k}^{i})}{c}\\
&= \frac{\tilde{v}_{ci}((t_{k+1}^{i})^{-})}{c}\geq\frac{\sqrt{\delta}}{c}.
\end{split}
\end{equation}

That is to say the distributed output-based event-triggered  mechanism \eqref{trigger2} has a minimal inter-execution time $\tau_{d}=\frac{\sqrt{\delta}}{c}$ which is strictly positive and independent of $k$. Thus, $\lim_{k \rightarrow \infty} t_{k}^{i} = \infty$, which
implies $T_{M}=\infty$. As a result,  the solution $\bar{x}_{c}(t)$ of $\eqref{xc3}$ exists for all time, which means that
the Zeno behavior cannot happen.

Next, choose $\delta=\frac{\lambda\lambda_{m}}{\lambda_{0}}\epsilon^{2}$. Then, according to \eqref{bxct2}, we have
\begin{equation}\label{xct2}
\begin{split}
\lim_{t\rightarrow+\infty}\sup \|\bar{x}_{c}(t)\|&\leq\!\lim_{t\rightarrow+\infty}\sup \sqrt{a\mathbf{e}^{-\alpha t} +b }\\
&=\sqrt{b }=\sqrt{\frac{\lambda_{0}\delta}{\lambda\lambda_{m}}}=\epsilon.\\
\end{split}
\end{equation}
Thus the distributed output feedback control law \eqref{bui3} together with the distributed output-based event-triggered mechanism \eqref{trigger2} solves CRPSP for the augmented system \eqref{system4a}.

It follows from Proposition \ref{Proposition1} that our proof is completed.
\end{Proof}


\begin{Remark}\label{RemarkTheorem1b}
Compared with some existing results on event-triggered cooperative output regulation problems in \cite{HuW2,Wang4,Wang5}, our approach has the following two features. First, since the control approaches used in \cite{HuW2,Wang4,Wang5} cannot handle the unknown parameters, they do not apply to our problem.
Second, the control laws in \cite{HuW2,Wang4,Wang5} are all piecewise continuous or continuous and thus cannot be directly implemented in digital platforms, while here our control law is piecewise constant and thus can be directly implemented in digital platforms.
\end{Remark}




\section{Simulations}\label{Example}
In this section, we provide an example to illustrate our design.
Consider a group of linear multi-agent systems of the form \eqref{system1}
with $z_{i}\in\mathbb{R}$, $\xi_{i}=\mbox{col}(\xi_{1i},\xi_{2i})\in\mathbb{R}^{2}$,
$[A_{1i}(w),A_{2i}(w),A_{3i}(w)]=[-8,2,3]+[w_{1i},w_{2i},w_{3i}]$, $[c_{1i}(w),c_{2i}(w),b_{i}(w)]=[-25,-10,2]+[w_{4i},w_{5i},w_{6i}]$, $E_{0i}(w)=[1,0]+[w_{7i},0],~E_{2i}(w)=[0,1]+[0,w_{8i}]$
for $i=1,2,3,4$.
 The exosystem takes the form of \eqref{exosystem1} with $S=\begin{bmatrix}
        \begin{smallmatrix}
            0 & 1 \\
            -1 & 0
        \end{smallmatrix}
    \end{bmatrix}$ and the output of the exosystem is given by $y_{0}=v_{1}$.
 It is easy to verify that Assumptions \ref{Ass1}-\ref{Ass2} are both satisfied. Here we assume that $|w_{ji}|\leq1$ for $i=1,2,3,4$ and $j=1,2\cdots,8$ and $|v_{j}|\leq1$ for $j=1,2$.

\begin{figure}[H]
\centering
\includegraphics[scale=0.5]{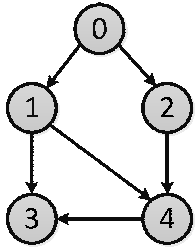}
\caption{Communication graph $\bar{\mathcal{G}}$} \label{graph}
\end{figure}

Figure \ref{graph} describes the  digraph $\bar{\mathcal{G}}$. Clearly, it contains a directed spanning tree with node 0 as the root. Thus, Assumption \ref{Ass3} is satisfied.
Note that  the minimal polynomial of $S$ is $s^2 + 1$. Then  we have $\Phi=\begin{bmatrix}
        \begin{smallmatrix}
            0 & 1 \\
            -1 & 0
        \end{smallmatrix}
    \end{bmatrix}$, $\Gamma=[1,0]$. 
Choose the controllable pair $(M,Q)$ with
$M=\begin{bmatrix}
        \begin{smallmatrix}
            0 & 1 \\
            -16 & -8
        \end{smallmatrix}
    \end{bmatrix},~Q=\begin{bmatrix}
        \begin{smallmatrix}
            0  \\
            1
        \end{smallmatrix}
    \end{bmatrix}$.
Solving the Sylvester equation \eqref{Sylvester} gives $\Gamma T^{-1}=[15, 8]$. By Theorem \ref{Theorem1}, we can obtain a distributed output feedback control law of the form \eqref{ui4} and a distributed output-based event-triggered mechanism \eqref{trigger2}, where various design parameters are
$\gamma_{0}=\gamma_{1}=1$, $\delta_{0}=4$, $\delta_{1}=4$,     $K=2$, $h=5$,
$A_{o}(h)=\begin{bmatrix}
        \begin{smallmatrix}
            -h\delta_{1} & 1 \\
            -h^{2}\delta_{0} & 0
        \end{smallmatrix}
    \end{bmatrix},~B_{o}(h)=\begin{bmatrix}
        \begin{smallmatrix}
            h\delta_{1} \\
            h^{2}\delta_{0}
        \end{smallmatrix}
    \end{bmatrix}$,
 and the design parameters $\bar{\sigma},\alpha$, $\beta$ and $\delta$ are given in Table \ref{table}.

\begin{figure}[H]
\centering
\includegraphics[scale=0.5]{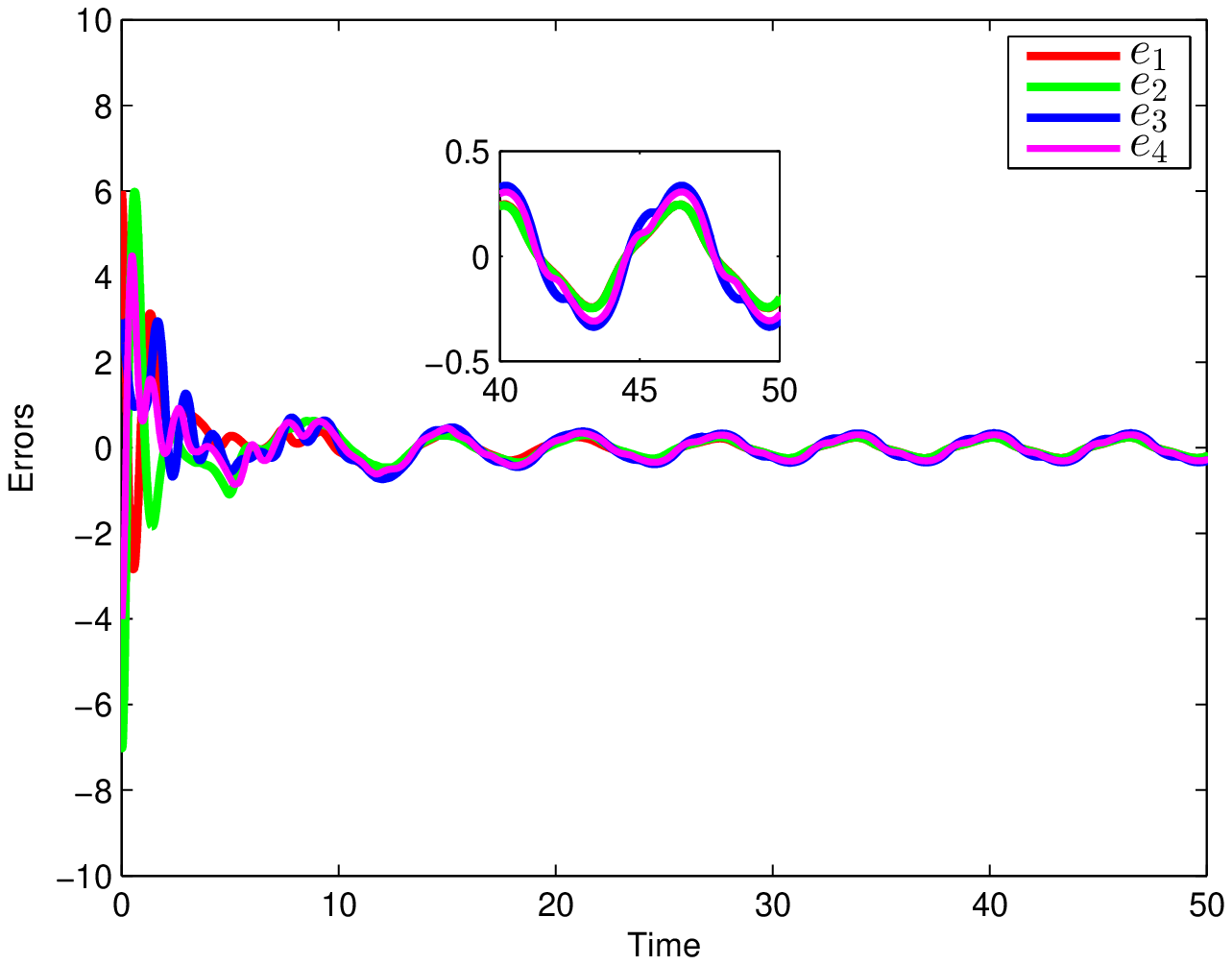}
\caption{Tracking errors of all followers for $t\in[0,50]$,} $\bar{\sigma}=0.06$, $\beta=3$, $\alpha=0.4$ and $\delta=0.001.$ \label{error1}
\end{figure}
\begin{figure}[H]
\centering
\includegraphics[scale=0.5]{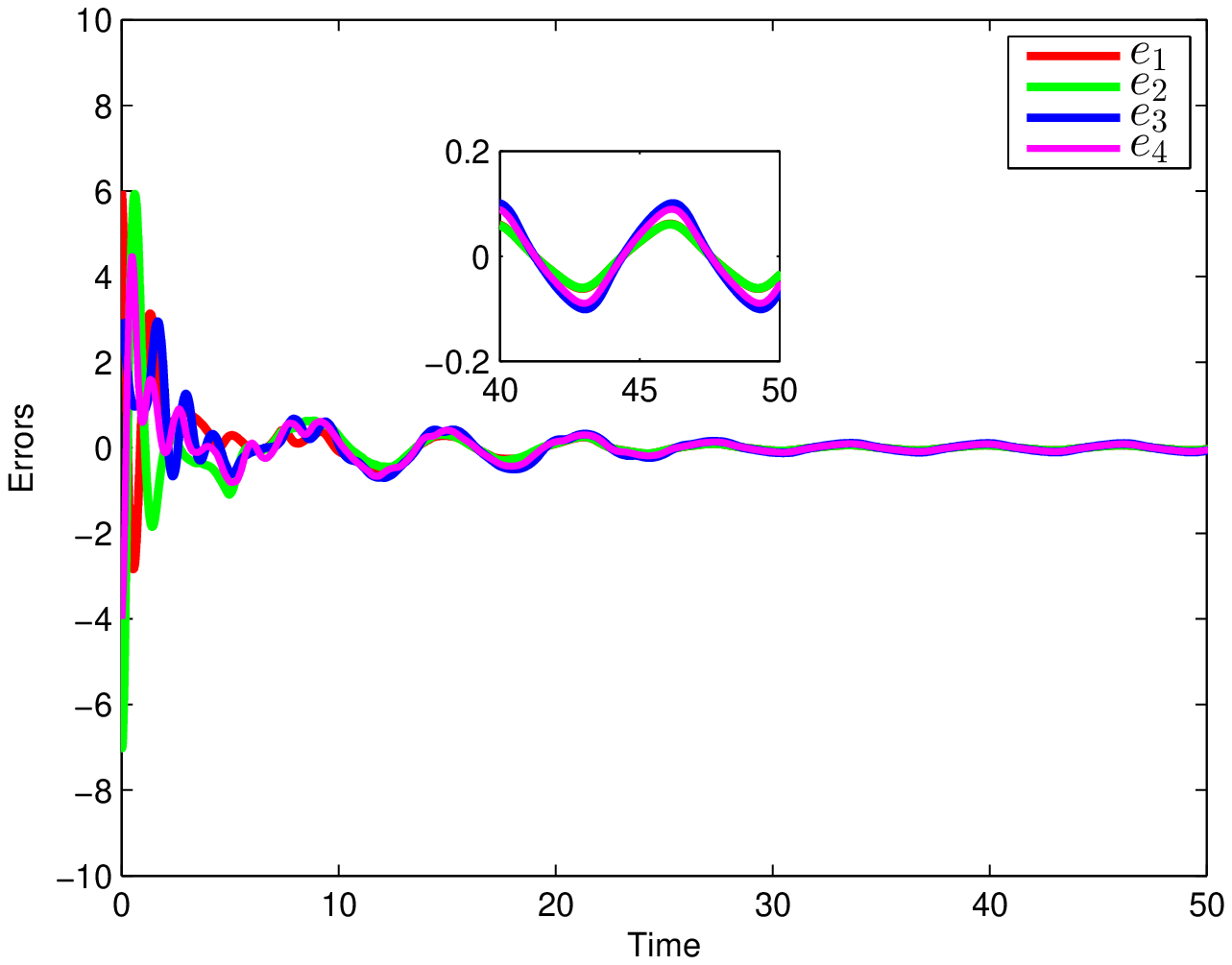}
\caption{Tracking errors of all followers for $t\in[0,50]$,} $\bar{\sigma}=0.06$, $\beta=3$, $\alpha=0.4$ and $\delta=0.0001.$ \label{error2}
\end{figure}
\begin{table}[H]
  \begin{center}
    \caption{Event-triggered numbers of all agents.}\label{table}
  \scalebox{0.78}{
    \begin{tabular}{|c|c|c|c|c|c|c|} \hline
      \multirow{2}{*}{Design parameters} & \multirow{2}{*}{Time} & \multicolumn{4}{c|}{Triggering numbers for each agent} \\
      \cline{3-6} &  & Agent 1 & Agent 2 & Agent 3 & Agent 4\\ \hline
      $\bar{\sigma}=0.06$, $\alpha=0.4$ & \multirow{2}{*}{0-5s}   & \multirow{2}{*}{59} & \multirow{2}{*}{71} & \multirow{2}{*}{104} & \multirow{2}{*}{110} \\
      $\beta=1$, $\delta=0.0001$ &    &  &  &  &  \\ \hline
      $\bar{\sigma}=0.06$, $\alpha=0.4$ & \multirow{2}{*}{0-5s}   & \multirow{2}{*}{50} & \multirow{2}{*}{59} & \multirow{2}{*}{86} & \multirow{2}{*}{84} \\
      $\beta=3$, $\delta=0.0001$ &    &  &  &  &  \\ \hline
      $\bar{\sigma}=0.06$, $\alpha=0.4$ & \multirow{2}{*}{45-50s}   & \multirow{2}{*}{128} & \multirow{2}{*}{133} & \multirow{2}{*}{176} & \multirow{2}{*}{137} \\
      $\beta=3$, $\delta=0.0001$ &    &  &  &  &  \\ \hline
      $\bar{\sigma}=0.06$, $\alpha=0.4$ & \multirow{2}{*}{45-50s}   & \multirow{2}{*}{38} & \multirow{2}{*}{40} & \multirow{2}{*}{72} & \multirow{2}{*}{64} \\
      $\beta=3$, $\delta=0.001$ &    &  &  &  &  \\ \hline
    \end{tabular}}
  \end{center}
\end{table}

 Simulation is performed with the unknown parameters
$w_{1}=[0.4,-0.3,0.6,-0.7$, $0.5,0.8,-0.2,-0.4]^{T}$, $w_{2}=[-0.3,0.2,0.5,-0.8,0.4,0.7$, $-0.1,-0.3]^{T}$, $w_{3}=[-0.6,-0.4,0.7,-0.8,0.4,0.2,-0.5$, $0.1]^{T}$, $w_{4}=[0.3$, $-0.4,0.5,-0.6,0.3,0.7,-0.4,0.2]^{T}$
and the initial conditions $z_{1}(0)=2,~z_{2}(0)=4,~z_{3}(0)=2,~z_{4}(0)=1$, $\xi_{1}(0)=[7,-5]^{T},~\xi_{2}(0)=[-6,-3]^{T}$, $\xi_{3}(0)=[4,-3]^{T},~\xi_{4}(0)=[-3,4]^{T},v(0)=[1,0]^{T}$, $\eta_{1}(0)=[ -3, 1  ]^{T},~\eta_{2}(0)=[3, -2]^{T}$,
$\eta_{3}(0)=[ 2, -1 ]^{T},~\eta_{4}(0)=[2, -1]^{T}$,
$\zeta_{1}(0)=[   2, -5]^{T},~\zeta_{2}(0)=[  2, 3]^{T}$,
$\zeta_{3}(0)=[  3, -4]^{T},~\zeta_{4}(0)=[2, -5]^{T}$.

Figures \ref{error1} and \ref{error2} show that, for $i=1,2,3,4$,  the tracking errors of all followers satisfy $\lim_{t\rightarrow\infty}|e_{i}(t)|\leq0.5$ for $\delta=0.001$ and $\lim_{t\rightarrow\infty}|e_{i}(t)|\leq0.2$ for $\delta=0.0001$. It can be seen that the larger $\delta$ leads to the larger tracking errors during the steady stage.
 Table \ref{table} shows the event-triggered numbers of all followers. It can be seen that the larger $\beta$ leads to the less triggering numbers during the initial stage, and the larger $\delta$ leads to the less triggering numbers during the steady stage.

%

\section{Conclusion}\label{Conclusion}
In this paper, 
we have studied the cooperative robust practical output regulation problem for a class of linear minimum-phase multi-agent systems. By designing a distributed  output feedback control law and a distributed output-based event-triggered mechanism, we have solved our problem and proved that the Zeno behavior can be prevented even if the neighbor-based measurement output is equal to zero.
As we have revealed in Remark \ref{RemarkProb1b} that
our event-triggered mechanism  relies on the continuous-time output information from neighboring agents, thus costing communication resources.
It is interesting to  further consider how to design novel triggering mechanisms that avoid continuous communication by appropriate samplings.

\section*{Acknowledgment}
The authors would like to thank the associate editor and the reviewers for their valuable comments and constructive suggestions.

\end{document}